 \documentclass[final,authoryear,1p,times]{elsarticle}
\usepackage{epsfig}

\usepackage{amssymb}
\usepackage{amsthm}
\usepackage{amsmath}

\usepackage{graphicx}
\usepackage{subfigure}
\usepackage{epstopdf}
\usepackage{color}

\usepackage{stfloats}
\usepackage{bbm}

\journal{}

\begin{document}

\begin{frontmatter}



\title{Curve and surface construction with moving B-splines}


\author{Xunnian Yang\corref{cor1}}
\date{}
%\ead[cor1]{Tel.: +86 571 87951609; fax: +86 571 87951428.}
\ead{yxn@zju.edu.cn}

\address{School of mathematical sciences, Zhejiang University, Hangzhou 310058, China}

\begin{abstract}
This paper proposes a simple technique of curve and surface construction with B-splines. Given a control polygon or a control mesh together with node ordinates corresponding to all control points, a rational curve or surface is obtained by least squares fitting of a moving constant to the control points with weights given by uniform B-splines centered at the prescribed nodes. This kind of curves and surfaces are natural generalizations of uniform B-spline curves and surfaces. By choosing proper nodes, the obtained curves can have sharp or rounded corners, partial or full straight edges while the obtained surfaces can have sharp or rounded vertices, sharp or smoothed edges, feature lines, etc. Except at sharp corners or sharp edges, the curves or surfaces have the same continuity orders as the moving B-splines. Practical examples have been given to demonstrate the effectiveness of the proposed technique for curve and surface modeling.
\end{abstract}


\end{frontmatter}

% \linenumbers


%%%%%%%%%%%%%%%%%%%%%%%%%%%%%%%%%%%%%%%%%%%%%%%%%%%%%%%%%%%%%%%%%%%%%%%%%%%%%%%%%%%%%%%%%%%
%%% Sectioin 1                                                                          %%%
%%%%%%%%%%%%%%%%%%%%%%%%%%%%%%%%%%%%%%%%%%%%%%%%%%%%%%%%%%%%%%%%%%%%%%%%%%%%%%%%%%%%%%%%%%%

\section{Introduction}
\label{Sec:intro}

Due to their nice properties such as local adjustability, shape-preserving, etc., B-spline curves and surfaces are prevalent tools for free-form shape modeling~\citep{Farin2001CAGDBook}. As the rational extensions of B-spline curves and surfaces, non-uniform rational B-spline (NURBS) curves and surfaces can represent analytical curves and surfaces like conics or quadrics exactly~\citep{Farin1999NURBSFromProjGeometry}. The ability of representing free-form curves and surfaces as well as analytical curves and surfaces in a unified way makes NURBS de facto industry standard in CAD/CAM commercial software~\citep{Piegl1997TheNURBSBook}. In addition to interactive design, NURBS curves and surfaces have also been successfully applied in many other fields such as reverse engineering~\citep{MaKruth1998NURBSfitting} or isogeometric analysis~\citep{Hughes2005:IsogeometricAnalysis}, etc.

For specific modeling purposes or even more powerful modeling results, techniques of B-spline or NURBS curves and surfaces have been extended much in the past few decades. Non-polynomial B-splines that can represent helices or conics exactly without rational form have been proposed in~\citep{PottmannW94:helixsplines,ZhangJiwen96:C-curves,LuWY02:uniformHyperbolicBspline}. By introducing additional degrees of freedom or shape parameters into the basis functions, the fullness or the proximity to the control polygon can be controlled efficiently~\citep{KovacsV17:P-curvesandSurfaces,KovacsV18:P-BezierandP-Bspline,Roth2021:CyclicProxityCS}. Recently, Yang proposed to replace the conventional real weights of NURBS curves and surfaces by matrices such that the obtained curves and surfaces can be controlled or edited using control normals or control tangents in addition to the control points~\citep{Yang16MatrixRational,Yang2018:fairingMatNURBScurve}. Another recent generalization of NURBS was given by Taheri, et al. who decoupled the unique weight for each control point into independent weights for separate coordinates~\citep{Suresh2019:GNURBSdecouplingweights}. Actually, NURBS curves and surfaces with decoupled weights are just matrix weighted NURBS curves and surfaces of which the weights are diagonal matrices.

To break through the restrictions of aligned knots for tensor product B-spline or NURBS surfaces, Sederberg, et al. proposed point-based surfaces, and in particular, T-spline surfaces~\citep{Sederberg2003:Tspline}. The control meshes of T-spline surfaces can have T-junctions, then the knot insertion and the connection of two T-spline surfaces can be locally implemented. Wang et al. proposed the so called generalized NURBS surfaces by constructing the point-based surfaces with arbitrary parameter domain~\citep{WangQing2004:GNURBS}. Similar technique has been applied in~\citep{RunionsS11Partitionofunityparametrics} who introduced local shape parameters for every weighted basis for meta-modeling. Rational Gaussian curves and surfaces presented in~\citep{Goshtasby95:RationalGaussianCurvesSurfaces,Goshtasby03:RaGaussCircleSphere} are also point-based curves and surfaces. Differently from point-based B-spline curves and surfaces, rational Gaussian curves and surfaces have infinite order of smoothness. When sharp features are to be modeled on a smooth surface, one can employ the technique of curved knots~\citep{ShiYSP2014:CurveKnotSplineSurface} or by non-uniform subdivision surfaces ~\citep{Sederberg1998:NonuniformRecursiveSubdivision,LiXin16:Non-uniformCCsurfaces,KosinkaSD14:semi-sharpCreasesSubdivision}. A different approach to model salient or sharp features on a smooth surface is to deform the original surface with displacements controlled by vectors and splines~\citep{KosinkaSD15:ControlVectorsForSplines}.

In addition to spline curves and surfaces with equal degrees in the whole parameter domain, curves and surfaces can also be constructed by variable degree or multi-degree (MD) splines~\citep{Costantini2000:VariableDegreePolynomialsplines,SederbergZS03:MDsplines}. MD-splines have many similar properties with traditional B-splines, and they are potentially powerful tools for geometric modeling. Particularly, curves or surfaces with sharp corners, local straight segments or local flat patches, etc., can be represented more compactly than by conventional B-splines. Recently, independent algorithms have been developed for constructing multi degree bases or evaluating points and derivatives on the obtained MD-spline curves and surfaces~\citep{ShenW10:basisofMDsplines,LiHL12:geometricMDspline,BeccariCM17:onMDsplines,BeccariC19:CoxdeBoorrecurrenceforMDspline}.

In this paper we present a novel generalization of uniform B-spline curves and surfaces to rational curves and surfaces. Our approach is motivated by curves and surfaces constructed by moving least squares fitting~\citep{Levin98:appr.powerMLS}. This is based on the observation that the basis functions for a uniform B-spline curve are obtained by moving an initial B-spline function along the parameter domain and a uniform B-spline curve is just the solution to least squares fitting of a moving constant to the control points with weights chosen by the moving B-splines. We generalize this basic procedure by using arbitrary chosen nodes for the moving least squares fitting. As a result, a moving B-spline-based rational curve is constructed from a sequence of given control points and the prescribed nodes corresponding to the control points. The same technique can be applied for surface construction with moving B-splines. Curves constructed by moving B-splines can have sharp or rounded corners, local straight edges, etc., while surfaces constructed by moving B-splines can have various features like sharp corners, feature lines, etc. just by choosing proper nodes for the control points. Similar to curves and surfaces constructed with MD-splines, curves and surfaces constructed with moving B-splines can represent complex shapes more compactly than by conventional NURBS curves and surfaces. Differently from MD-splines, curve and surface construction with moving B-splines benefits the advantages of simplicity and efficiency.

The remaining of the paper is structured as follows. In Section~\ref{Sec:Moving_Bspline_curve} we propose techniques of curve construction with moving B-splines. Properties and techniques for feature modeling on a moving B-spline curve will also be discussed or given. Section~\ref{Sec:Moving_Bspline_surface} is devoted to the construction of surfaces with moving B-splines. Section~\ref{Sec:Examples} presents several interesting examples for curve and surface modeling by the proposed technique and the paper is concluded with a brief summary in Section~\ref{Sec:Conclude}.

%%%%%%%%%%%%%%%%%%%%%%%%%%%%%%%%%%%%%%%%%%%%%%%%%%%%%%%%%%%%%%%%%%%%%%%%%%%%%%%%%%%%%%%%%%%
%%% Sectioin 2                                                                          %%%
%%%%%%%%%%%%%%%%%%%%%%%%%%%%%%%%%%%%%%%%%%%%%%%%%%%%%%%%%%%%%%%%%%%%%%%%%%%%%%%%%%%%%%%%%%%

\section{Curve construction with moving B-splines}
\label{Sec:Moving_Bspline_curve}

This section proposes techniques of curve modeling with moving B-splines. Properties of such curves will be discussed and presented.

\subsection{Curve modeling with moving B-splines}
Assume $\tau_0\leq \tau_1\leq \ldots\leq \tau_k$ and $\tau_0\neq \tau_k$. A B-spline of order $k$ (degree $k-1$), namely $N_{0,k}(t)$, with support $[\tau_0,\tau_k]$, is defined recursively by
\begin{eqnarray*}
  N_{i,1}(t) &=& \left\{
                    \begin{array}{ll}
                      1, & \hbox{if $t\in[\tau_i,\tau_{i+1})$} \\
                      0, & \hbox{otherwise}
                    \end{array}
                  \right.
   \\
  N_{i,k}(t) &=& {{t-\tau_i}\over{\tau_{i+k-1}-\tau_i}}N_{i,k-1}(t) + {{\tau_{i+k}-t}\over{\tau_{i+k}-\tau_{i+1}}}N_{i+1,k-1}(t)
\end{eqnarray*}
with $0/0=0$ if $\tau_i=\tau_{i+1}$.

To construct curves with even fewer parameters, in this paper we assume that the original B-splines are defined by uniform knots. In particular, a B-spline of order $k$ is defined by choosing the knots as $\tau_0=-\frac{k}{2}, \tau_1=-\frac{k}{2}+1, \ldots, \tau_k=\frac{k}{2}$. Then the obtained B-spline $N_{0,k}(t)$ has a symmetric support $[-\frac{k}{2},\frac{k}{2}]$. For simplicity purpose, $N_{0,k}(t)$ will be denoted as $N_{k}(t)$ in the following text.

\subsubsection{Moving B-spline curves}
Suppose that $P_0, P_1, \ldots, P_n$ are a sequence of points lying in $\mathbb{R}^d$ and $t_i\in\mathbb{R}$, $i=0,1,\ldots,n$, are the given nodes corresponding to the points. We construct a curve by least squares fitting of a moving constant to the points with weights given by the moving B-splines centered at the prescribed nodes. Let
\[
F(P(t))=\sum_{i=0}^n N_k(t-t_i)[P(t)-P_i]^2.
\]
By minimizing $F(P(t))$ or solving the equation $\frac{\partial F(P(t))}{\partial P(t)}=0$, we have
\begin{equation}\label{Eqn:Moving_Bspline_curve}
P(t)=\frac{\sum_{i=0}^n P_i N_k(t-t_i)}{\sum_{i=0}^n N_k(t-t_i)}.
\end{equation}
Though the nodes can be arbitrary real numbers, we assume further that the nodes satisfy $t_0<t_1<\ldots<t_n$ and $t_{i+1}-t_i<k$, $i=0,1,\ldots,n-1$, such that the curve $P(t)$ is a continuous curve on the domain $[t_1-\frac{k}{2},t_{n-1}+\frac{k}{2}]$.

Let
\[
R_i(t)=\frac{N_k(t-t_i)}{\sum_{j=0}^n N_k(t-t_j)}, \ \ \ i=0,1,\ldots,n.
\]
It yields that $R_i(t)\geq 0$, $i=0,1,\ldots,n$, and $\sum_{i=0}^n R_i(t)\equiv1$. It is also easily verified that the functions $R_i(t)$ are linear independent over the domain. Now, Equation~(\ref{Eqn:Moving_Bspline_curve}) can be reformulated as
\[
P(t)=\sum_{i=0}^n P_i R_i(t), \ \ \ t\in\left[t_1-\frac{k}{2},t_{n-1}+\frac{k}{2}\right].
\]
In particular, if the knots satisfy $t_{i+1}=t_i+1$, $i=0,1,\ldots,n-1$, the curve $P(t)$ will reduce to a uniform B-spline curve on the domain $[t_0+\frac{k}{2}-1,t_n-\frac{k}{2}+1]$. This implies that curves constructed with moving B-splines are natural generalizations of uniform B-spline curves. We also refer this kind of curves as \emph{moving B-spline curves}. Like conventional rational B-spline curves, moving B-spline curves also have nice properties like convex hull property, variation diminishing property, etc. The properties together with freely chosen nodes will make the moving B-spline curves a powerful tool for geometric modeling.

\begin{figure}[htb]
  \centering
  \subfigure[]{\includegraphics[width=0.45\linewidth]{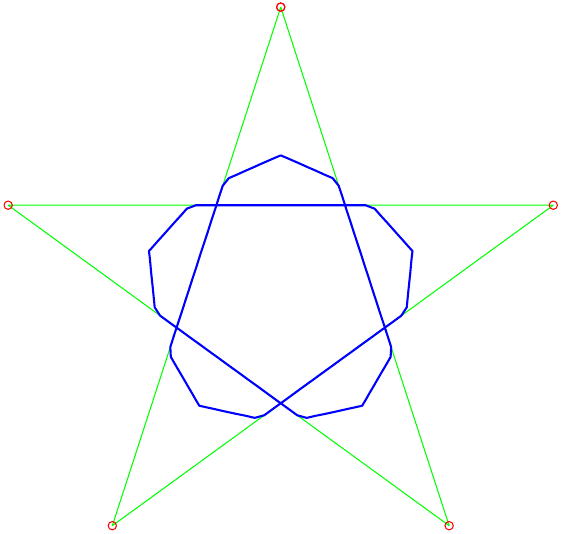}}\\
  \subfigure[]{\includegraphics[width=0.45\linewidth]{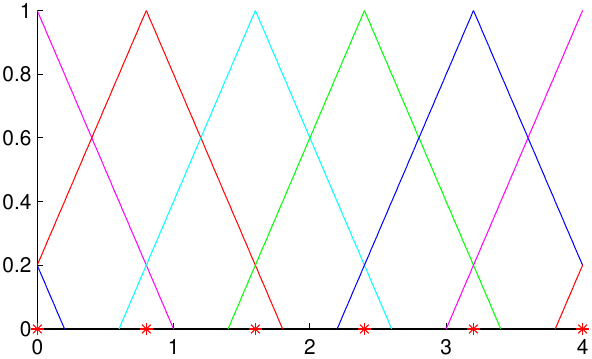}}
  \subfigure[]{\includegraphics[width=0.45\linewidth]{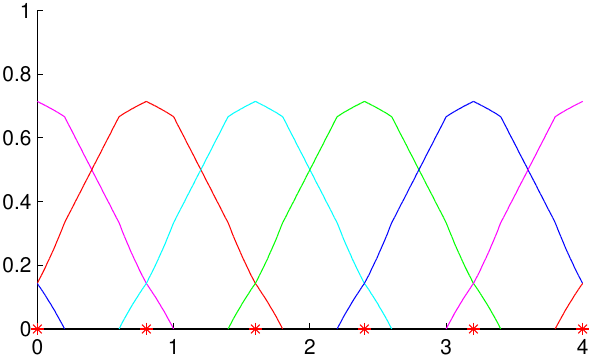}}
  \caption{(a) A closed moving B-spline curve of degree 1 together with control polygon; (b) the moving linear B-spline bases; (c) the normalized moving B-spline bases. The stars '\textcolor{red}{*}' denote the prescribed nodes for the moving B-splines.}
  \label{Fig:moving linear_Bspline curves_star5}
\end{figure}

For arbitrary $t\in[t_1-\frac{k}{2},t_{n-1}+\frac{k}{2}]$, $P(t)$ can be computed by
\begin{equation}\label{Eqn:Moving_Bspline_curve_point}
P(t)=\frac{\sum_{|t-t_i| <\frac{k}{2}} P_i N_k(t-t_i)}{\sum_{|t-t_i| <\frac{k}{2}} N_k(t-t_i)}.
\end{equation}
In particular, the two end points of the curve are obtained as
\[P(t_1-\frac{k}{2})=P_0\]
and
\[P(t_{n-1}+\frac{k}{2})=P_n.\]

\begin{figure}[htb]
  \centering
  \subfigure[]{\includegraphics[width=0.45\linewidth]{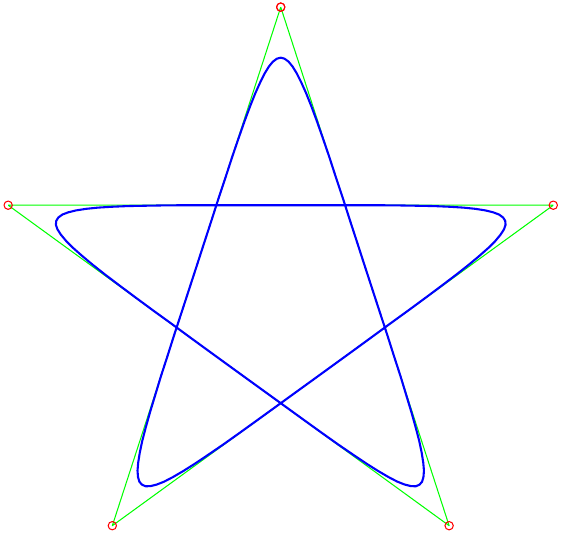}}\\
  \subfigure[]{\includegraphics[width=0.45\linewidth]{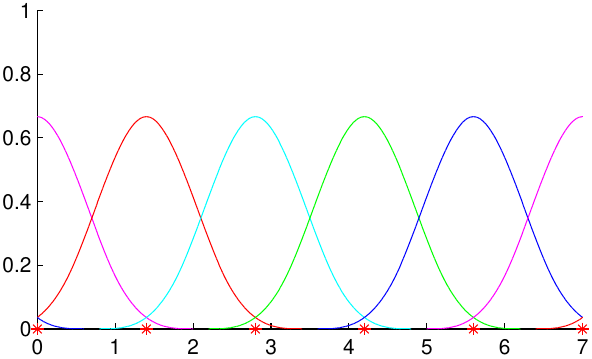}}
  \subfigure[]{\includegraphics[width=0.45\linewidth]{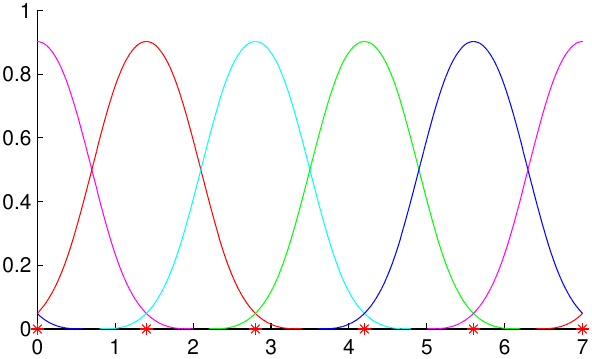}}
  \caption{(a) A closed moving B-spline curve of degree 3 together with control polygon; (b) the moving cubic B-spline bases; (c) the normalized moving B-spline bases.}
  \label{Fig:moving cubic_Bspline curves_star5}
\end{figure}

If the control points satisfy $P_n=P_0$, a closed curve will be constructed from the points with the moving B-spline. Assume $t_i$, $i=0,1,\ldots,n$, are the nodes corresponding to the control points, a closed moving B-spline curve is given by
\begin{equation}\label{Eqn:closed_Moving_Bspline_curve}
\begin{array}{lll}
P(t)&=&\frac{\sum_{i=0}^{n-1} P_i N_k(t-t_i)}{\sum_{i=0}^{n-1} N_k(t-t_i)}, \\
    &=&\sum_{i=0}^{n-1}P_i R_i(t), \ \ \ t\in[t_0,t_n],
\end{array}
\end{equation}
where $R_i(t)=\frac{N_k(t-t_i)}{\sum_{j=0}^{n-1} N_k(t-t_j)}$, $i=0,1,\ldots,n-1$.

Let $T=t_n-t_0$. For arbitrary $t\in[t_0,t_n]$, the point $P(t)$ is computed by
\begin{equation}\label{Eqn:closed_Moving_Bspline_curve_point}
P(t)=\frac{\sum_{\delta t_i<\frac{k}{2}} P_i N_k(t-t_i)}{\sum_{\delta t_i<\frac{k}{2}} N_k(t-t_i)},
\end{equation}
where $\delta t_i=\min\{|t-t_i|,|t-t_i+T|,|t-t_i-T|\}$.

Figure~\ref{Fig:moving linear_Bspline curves_star5} illustrates a closed curve constructed by moving B-splines of degree 1. Given a closed control polygon consisting of 5 points, the moving B-spline curve of degree 1 is constructed by choosing nodes as $t_i=t_{i-1}+0.8$, $i=1,2,\ldots,5$. Unlike conventional B-spline curves of degree 1 that coincide with the control polygons exactly, the normalized moving B-splines of degree 1 can be employed to smooth out the control polygon. Using the same set of control points as that illustrated in Figure~\ref{Fig:moving linear_Bspline curves_star5}, another closed curve is constructed by moving B-splines of degree 3. In particular, by choosing the nodes as $t_i=t_{i-1}+1.4$, $i=1,2,\ldots,5$, only two moving cubic B-splines do not vanish on the intervals $[t_i+2,t_{i+4}-2]$. As a result, a $C^2$ continuous cubic moving B-spline curve that has local straight line segments is obtained; see Figure~\ref{Fig:moving cubic_Bspline curves_star5}.

\subsubsection{Modeling curve features with moving B-splines}
By choosing appropriate nodes corresponding to the control points, the obtained moving B-spline curves can have various types of features. A control point can be a sharp corner or a rounded corner while an edge can be a full straight edge or a partial straight edge embedded in the curve.

$\bullet$ \textbf{Modeling sharp corners.}
Assume $P(t)$ is a moving B-spline curve of order $k$ given by Equation~(\ref{Eqn:Moving_Bspline_curve}) or Equation~(\ref{Eqn:closed_Moving_Bspline_curve}). If the nodes satisfy $t_{i-1}+\frac{k}{2} \leq t_{i+1}-\frac{k}{2}$, or equivalently, $t_{i+1}-t_{i-1}\geq k$, it is verified that $R_i(t)=1$ for $t\in[t_{i-1}+\frac{k}{2}, t_{i+1}-\frac{k}{2}]$ but $R_j(t)=0$, $j\neq i$, in the same interval. It thus follows that $P(t)=P_i$ when $t\in[t_{i-1}+\frac{k}{2}, t_{i+1}-\frac{k}{2}]$.

It is also easily verified that only two bases $R_{i-1}(t)\neq 0$ and $R_i(t)\neq 0$ when $t\in[t_i-\frac{k}{2},t_{i-1}+\frac{k}{2}]\bigcap[t_{i-2}+\frac{k}{2},t_{i-1}+\frac{k}{2}]$. It just implies that the curve $P(t)$ is local straight along the edge $P_{i-1}P_i$ when $t$ belongs to the mentioned parameter interval. Similarly, the curve $P(t)$ is local straight in the interval $[t_{i+1}-\frac{k}{2},t_i+\frac{k}{2}]\bigcap[t_{i+1}-\frac{k}{2},t_{i+2}-\frac{k}{2}]$. Combining these together, we know that $P_i$ is a sharp corner bounded by two local straight edges.

\begin{figure}[htb]
  \centering
  \subfigure[]{\includegraphics[width=0.4\linewidth]{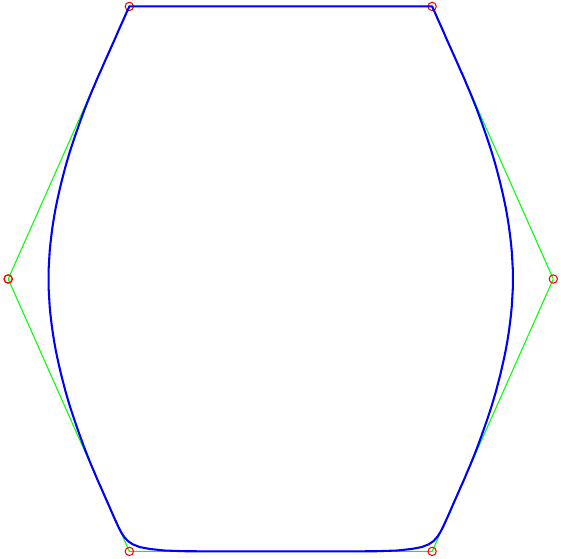}}\\
  \subfigure[]{\includegraphics[width=0.49\linewidth]{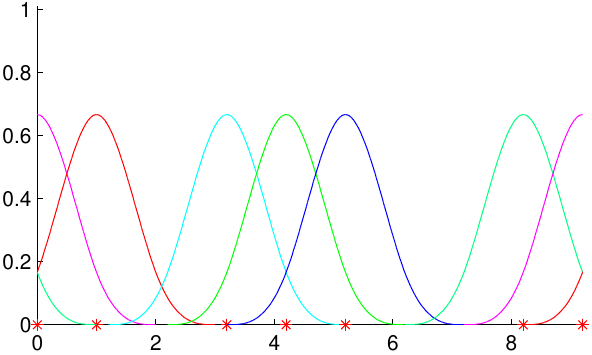}}
  \subfigure[]{\includegraphics[width=0.49\linewidth]{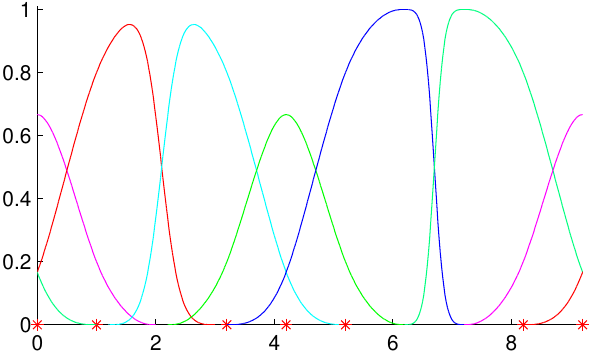}}
  \caption{A moving B-spline curve of degree 3 with non-uniform nodes: (a) the moving B-spline curve; (b) the moving cubic B-spline bases; (c) the normalized moving B-spline bases. }
  \label{Fig:moving Bspline curves_six}
\end{figure}

$\bullet$ \textbf{Modeling rounded corners.}
Contrary to the case of sharp corners, if the nodes of a moving B-spline curve of order $k$ satisfy $t_{i-1}+\frac{k}{2} > t_{i+1}-\frac{k}{2}$, or equivalently, $t_{i+1}-t_{i-1}< k$, the intersection of supports of $R_{i-1}(t)$, $R_i(t)$ and $R_{i+1}(t)$ is no longer null. In other words, at least three basis functions satisfy $R_{i-1}(t_i)>0$, $R_i(t_i)>0$ and $R_{i+1}(t_i)>0$. It follows that the point $P(t_i)$ is a weighted average of point $P_i$ and its neighbor points. In this case, we call $P_i$ a rounded corner of the curve.

$\bullet$ \textbf{Modeling full straight edges.}
If both of the two end points of an edge of the control polygon of a moving B-spline curve are sharp corners, the straight edge lies on the moving curve entirely. We refer the edge as the full straight edge of the curve.

Suppose $P(t)$ is a moving B-spline curve of order $k$. If the nodes satisfy $t_{i+1}-t_{i-1}\geq k$ and $t_{i+2}-t_i\geq k$, both points $P_i$ and $P_{i+1}$ are sharp corners. The curve segment $P(t)$, $t_i\leq t \leq t_{i+1}$, coincides with the edge $P_iP_{i+1}$ exactly.

$\bullet$ \textbf{Modeling partial straight edges.}
If there exists an interval, namely $[t_a,t_b]$, on which only two bases satisfy $R_i(t)>0$ and $R_{i+1}(t)>0$, then the curve segment $P(t)$, $t\in[t_a,t_b]$, is a straight line on the edge $P_iP_{i+1}$. This line segment is also referred as partial straight edge.

Practically, if the nodes satisfy $t_{i+1}-t_i<k$ and $t_{i+2}-t_{i-1}>k$, then the curve segment $P(t)$, $t\in[t_a,t_b]$, where $t_a=\max\{t_{i-1}+\frac{k}{2},t_{i+1}-\frac{k}{2}\}$ and $t_b=\max\{t_{i}+\frac{k}{2},t_{i+2}-\frac{k}{2}\}$, is locally straight. If $t_{i+1}-t_{i-1}<k$ or $t_{i+2}-t_i<k$, the point $P_i$ or $P_{i+1}$ is a rounded corner and the curve segment $P(t)$, $t\in[t_a,t_b]$ is strictly a partial straight edge.

Figure~\ref{Fig:moving Bspline curves_six} illustrates a closed curve constructed by moving cubic B-splines with control points $(-0.9,0)$, $(-0.5,-0.9)$, $(0.5,-0.9)$, $(0.9,0)$, $(0.5,0.9)$, $(-0.5,0.9)$. The prescribed nodes for the curve are $(0, 1, 3.2, 4.2, 5.3, 8.2, 9.2)$. From the figure we can see that the top two vertices are sharp corners while the remaining vertices are the rounded corners. Meanwhile, the top edge is a full straight edge and the bottom edge is a partial straight edge of the obtained curve due to different lengths of intervals between the neighboring nodes. This example also demonstrates that the moving B-spline curve can lie close or on an edge when the node interval has been increased.

\begin{figure}[htb]
  \centering
  \subfigure[]{\includegraphics[width=0.4\linewidth]{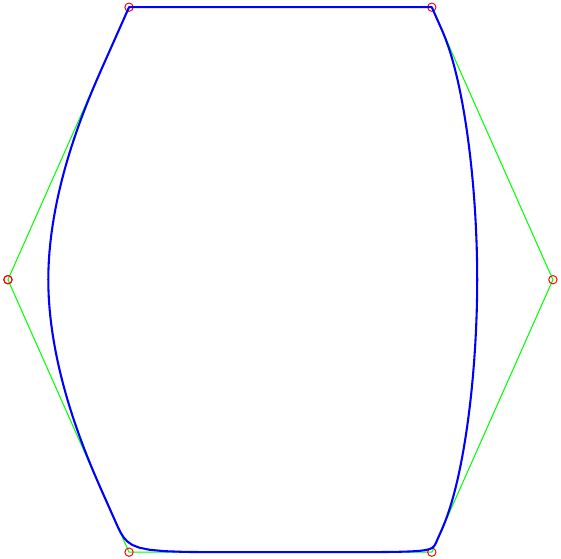}} \hskip 0.6cm
  \subfigure[]{\includegraphics[width=0.4\linewidth]{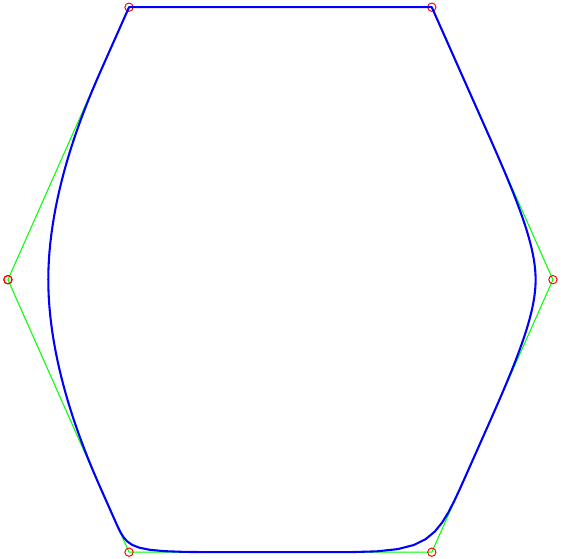}}
  \caption{(a) A weighted moving B-spline curve of degree 3 with weights $\omega_i=1$ but $\omega_3=0.3$; (b) A weighted moving B-spline curve of degree 3 with weights $\omega_i=1$ but $\omega_3=3$. }
  \label{Fig:moving Bspline curves_six_rational}
\end{figure}

\subsection{Weighted moving B-spline curves}
Similar to conventional rational B-spline curves, moving B-spline curves can have weights. Suppose that $P_0, P_1, \ldots, P_n$ are a sequence of given points. Assume $t_0<t_1<\ldots<t_n$ are the corresponding nodes and $\omega_i>0$, $i=0,1,\ldots,n$, are a set of given weights. A weighted moving B-spline curve of order $k$ is given by
\begin{equation}\label{Eqn:weighted_Moving_Bspline_curve}
P(t)=\frac{\sum_{i=0}^n \omega_i P_i N_k(t-t_i)}{\sum_{i=0}^n \omega_i N_k(t-t_i)}, \ \ \ t\in\left[t_1-\frac{k}{2},t_{n-1}+\frac{k}{2}\right].
\end{equation}
Let $R_i^\omega(t)=\frac{\omega_i N_k(t-t_i)}{\sum_{j=0}^n \omega_j N_k(t-t_j)}$, $i=0,1,\ldots,n$. The weighted moving B-spline curve can be reformulated as $P(t)=\sum_{i=0}^n P_i R_i^\omega(t)$. In particular, if a selected weight, namely $\omega_{i_0}$, approaches infinity, the function $R_i^\omega(t)$ approaches 1 for $t\in (t_i-\frac{k}{2},t_i+\frac{k}{2})$. Therefore, the curve $P(t)$ approaches the point $P_{i_0}$ when $t\in (t_i-\frac{k}{2},t_i+\frac{k}{2})$. This just implies that besides the control points and nodes, the weights can also be used for shape editing of moving B-spline curves.
Figure \ref{Fig:moving Bspline curves_six_rational} illustrates two weighted moving B-spline curves modified from the moving B-spline curve as defined in Figure \ref{Fig:moving Bspline curves_six}. By choosing all weights $\omega_i=1$ except for $\omega_3=0.3$ or $\omega_3=3$, the obtained curve will be locally pulled far away or closer to the corresponding control point.

%%%%%%%%%%%%%%%%%%%%%%%%%%%%%%%%%%%%%%%%%%%%%%%%%%%%%%%%%%%%%%%%%%%%%%%%%%%%%%%%%%%%%%%%%%%
%%% Sectioin 3                                                                          %%%
%%%%%%%%%%%%%%%%%%%%%%%%%%%%%%%%%%%%%%%%%%%%%%%%%%%%%%%%%%%%%%%%%%%%%%%%%%%%%%%%%%%%%%%%%%%

\section{Surface construction with moving B-splines}
\label{Sec:Moving_Bspline_surface}

This section presents techniques of surface modeling with moving bivariate B-splines. In particular, we are interested in surface modeling with moving tensor product B-splines. Salient features on the moving tensor product B-spline surfaces can be modeled using similar techniques as moving B-spline curves. The generalization of moving B-spline surfaces to weighted moving B-spline surfaces will also be given.

\subsection{Moving B-spline surfaces}
Let $k_1$ and $k_2$ are two positive integers. Assume that $N_{k_1}(s)$ is a uniform B-spline of order $k_1$ defined by the knot vector $\{-\frac{k_1}{2},-\frac{k_1}{2}+1,\ldots,\frac{k_1}{2}\}$
and $N_{k_2}(t)$ is a uniform B-spline of order $k_2$ defined by the knot vector $\{-\frac{k_2}{2},-\frac{k_2}{2}+1,\ldots,\frac{k_2}{2}\}$.

Suppose that $P_{ij}\in\mathbb{R}^d$, $i=0,1,\ldots,m$, $j=0,1,\ldots,n$, are a set of given points and $(s_{ij},t_{ij})$ are the nodes corresponding to the points. Similar to the definition of moving B-spline curves, a moving B-spline surface is given by
\begin{equation}\label{Eqn:Moving_Bspline_surface}
%\begin{array}{l}
P(s,t)=\frac{\sum_{i=0}^m\sum_{j=0}^n P_{ij} N_{k_1}(s-s_{ij})N_{k_2}(t-t_{ij})}{\sum_{i=0}^m\sum_{j=0}^n N_{k_1}(s-s_{ij})N_{k_2}(t-t_{ij})}.
%\end{array}
\end{equation}
Let $D_{ij}=[s_{ij}-\frac{k_1}{2},s_{ij}+\frac{k_1}{2}]\times[t_{ij}-\frac{k_2}{2},t_{ij}+\frac{k_2}{2}]$. The parameter domain of the surface~(\ref{Eqn:Moving_Bspline_surface}) is the union of supports of all moving B-splines, i.e., $D=\bigcup_{ij}D_{ij}$. In particular, if the nodes satisfy $s_{ij}=s_i$, $t_{ij}=t_j$, $i=0,1,\ldots,m$, $j=0,1\ldots,n$, the parameter domain of the surface~(\ref{Eqn:Moving_Bspline_surface}) will be
\[
D=\left[s_1-\frac{k_1}{2},s_{m-1}+\frac{k_1}{2}\right]\times\left[t_1-\frac{k_2}{2},t_{n-1}+\frac{k_2}{2}\right].
\]

For any given parameter pair $(s,t)$ within the domain, the point $P(s,t)$ can be computed by
\begin{equation}
P(s,t)=\frac{\sum_{|s-s_{ij}| <\frac{k_1}{2},|t-t_{ij}|<\frac{k_2}{2}} P_{ij} N_{k_1}(s-s_{ij})N_{k_2}(t-t_{ij})}{\sum_{|s-s_{ij}|<\frac{k_1}{2},|t-t_{ij}| <\frac{k_2}{2}} N_{k_1}(s-s_{ij})N_{k_2}(t-t_{ij})}.
\end{equation}
If the parameter domain is $D=[s_1-\frac{k_1}{2},s_{m-1}+\frac{k_1}{2}]\times[t_1-\frac{k_2}{2},t_{n-1}+\frac{k_2}{2}]$, the four corners of the surface patch are just the corner control points. It yields that
\[
\begin{array}{l}
P(s_1-\frac{k_1}{2},t_1-\frac{k_2}{2})=P_{00},  \\
P(s_{m-1}+\frac{k_1}{2},t_1-\frac{k_2}{2})=P_{m0}, \\
P(s_1-\frac{k_1}{2},t_{n-1}+\frac{k_2}{2})=P_{0n},  \\
P(s_{m-1}+\frac{k_1}{2},t_{n-1}+\frac{k_2}{2})=P_{mn}.
\end{array}
\]

If the control points satisfy $P_{0j}=P_{mj}$, $s_{mj}-s_{0j}=s_{m0}-s_{00}$, $j=0,1,\ldots,n$, the moving B-spline surface is obtained as follows:
\begin{equation}\label{Eqn:Moving_Bspline_surface_close_u}
\begin{array}{l}
P(s,t)=\frac{\sum_{i=0}^{m-1}\sum_{j=0}^n P_{ij} N_{k_1}(s-s_{ij})N_{k_2}(t-t_{ij})}{\sum_{i=0}^{m-1}\sum_{j=0}^n N_{k_1}(s-s_{ij})N_{k_2}(t-t_{ij})}.
%(s,t)\in[s_0,s_m]\times[t_1-\frac{k_2}{2},t_{n-1}+\frac{k_2}{2}].
\end{array}
\end{equation}
This surface is closed in $s$ direction.
If the control points satisfy $P_{i0}=P_{in}$, $t_{in}-t_{i0}=t_{0n}-t_{00}$, $i=0,1,\ldots,m$, the obtained moving B-spline surface becomes
\begin{equation}\label{Eqn:Moving_Bspline_surface_close_v}
\begin{array}{l}
P(s,t)=\frac{\sum_{i=0}^m\sum_{j=0}^{n-1} P_{ij} N_{k_1}(s-s_{ij})N_{k_2}(t-t_{ij})}{\sum_{i=0}^m\sum_{j=0}^{n-1} N_{k_1}(s-s_{ij})N_{k_2}(t-t_{ij})}.
%(s,t)\in[s_1-\frac{k_1}{2},s_{m-1}+\frac{k_1}{2}]\times[t_0,t_n].
\end{array}
\end{equation}
Now, the surface is closed in $t$ direction.
If the control points satisfy $P_{0j}=P_{mj}$, $s_{mj}-s_{0j}=s_{m0}-s_{00}$, $j=0,1,\ldots,n$, and $P_{i0}=P_{in}$, $t_{in}-t_{i0}=t_{0n}-t_{00}$, $i=0,1,\ldots,m$, a closed surface both in $s$ and $t$ directions is obtained as
\begin{equation}\label{Eqn:Moving_Bspline_surface_close_uv}
\begin{array}{l}
P(s,t)=\frac{\sum_{i=0}^{m-1}\sum_{j=0}^{n-1} P_{ij} N_{k_1}(s-s_{ij})N_{k_2}(t-t_j)}{\sum_{i=0}^{m-1}\sum_{j=0}^{n-1} N_{k_1}(s-s_{ij})N_{k_2}(t-t_{ij})}.
%(s,t)\in[s_0,s_m]\times[t_0,t_n].
\end{array}
\end{equation}
Evaluation of points on a closed surface in one or two directions can be employed the similar techniques for evaluating closed moving B-spline curves or evaluating open moving B-spline surfaces.

\begin{figure*}[htb]
  \centering
  \subfigure[]{\includegraphics[width=0.4\linewidth]{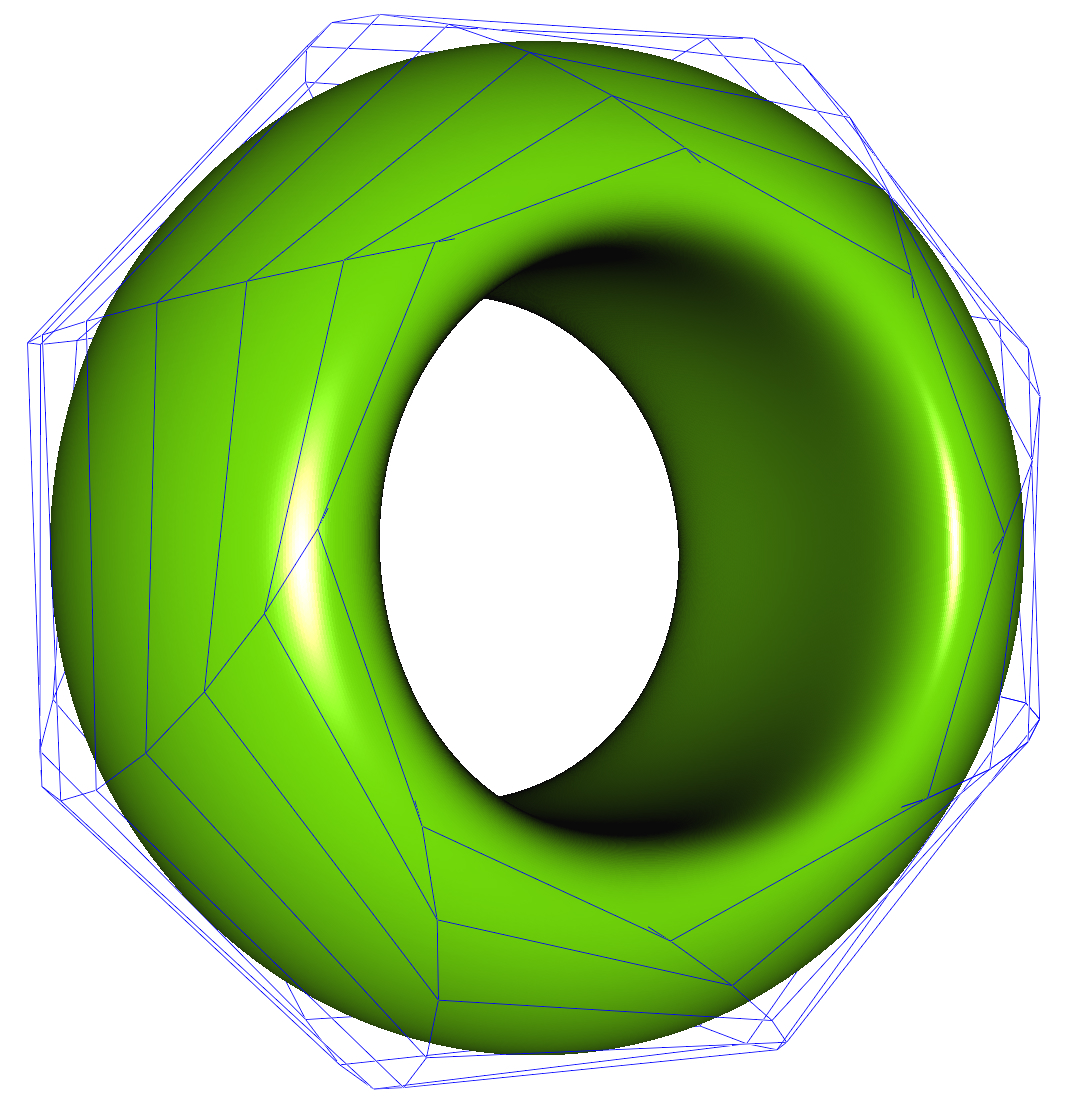}}\\
  \subfigure[]{\includegraphics[width=0.4\linewidth]{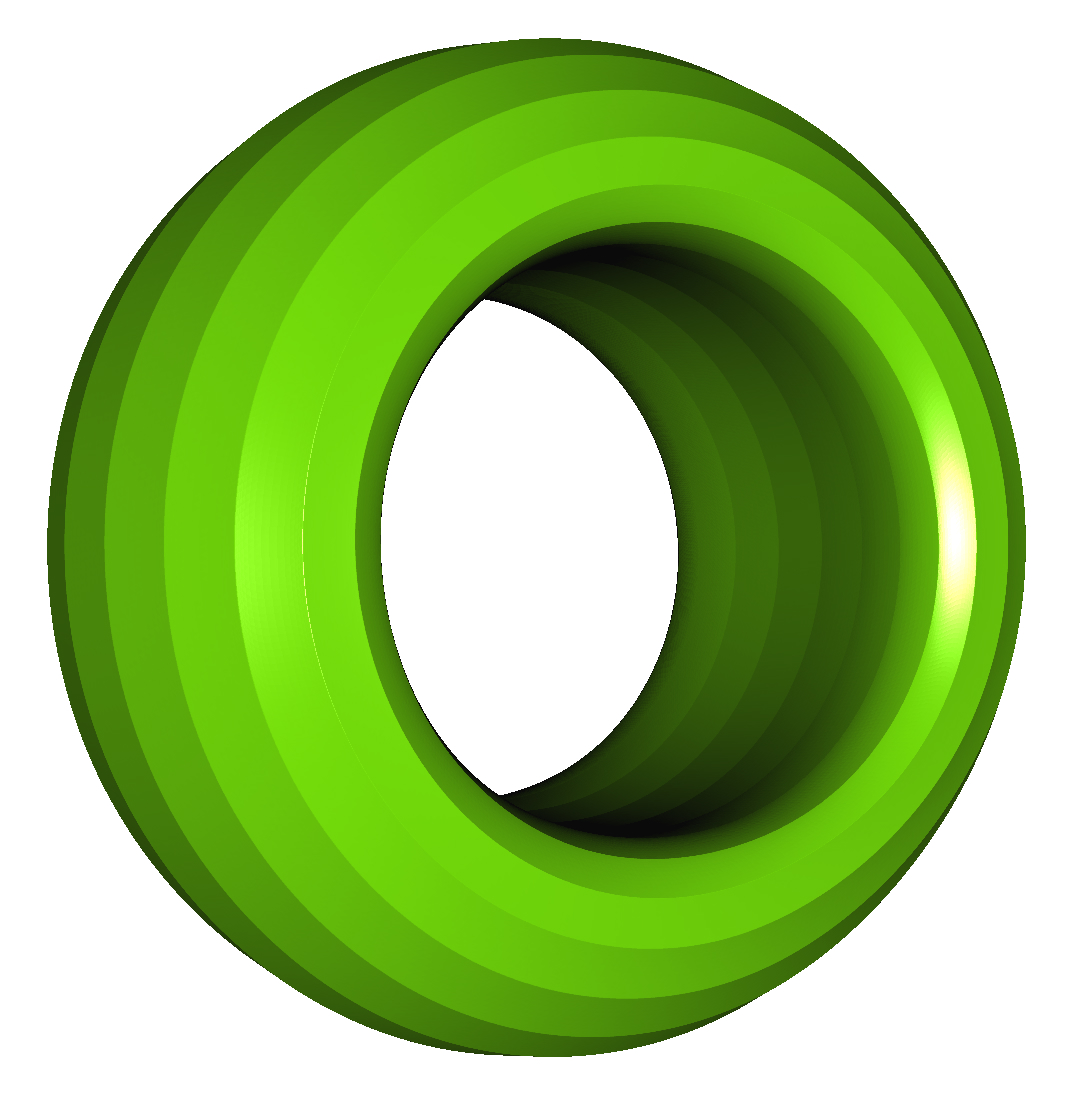}}
  \subfigure[]{\includegraphics[width=0.4\linewidth]{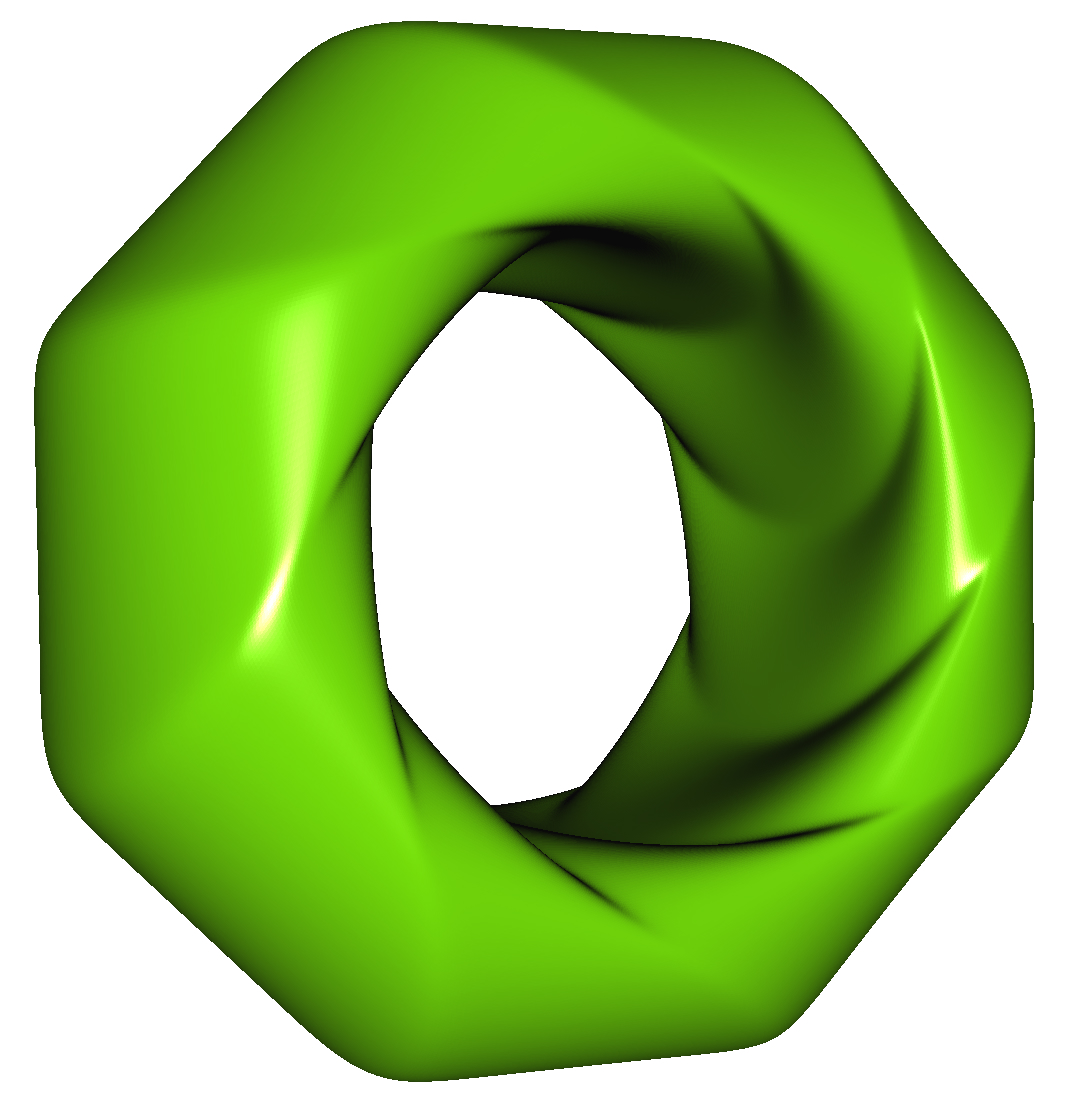}}
  \caption{(a) a bicubic B-spline surface consisting of $8\times20$ control points; (b) moving B-spline surface with nodes $s_i=s_{i-1}+1.0$ and $t_j=t_{j-1}+1.8$; (c) moving B-spline surface with nodes $s_i=s_{i-1}+1.5$ and $t_j=t_{j-1}+1.0$. }
  \label{Fig:moving Bspline surface_Ring}
\end{figure*}

Figure~\ref{Fig:moving Bspline surface_Ring} illustrates an example of moving B-spline surface construction from a given control mesh. A uniform bicubic B-spline surface with control mesh consisting of $8\times20$ control points has been given in Figure~\ref{Fig:moving Bspline surface_Ring}(a). When we choose the node vectors $s_{ij}=s_i=s_{i-1}+1.0$, $i=1,\ldots,8$, and $t_{ij}=t_j=t_{j-1}+1.8$, $j=1,\ldots,20$, a piecewise $C^2$ continuous surface with salient features has been obtained by employing Equation~(\ref{Eqn:Moving_Bspline_surface_close_uv}); see Figure~\ref{Fig:moving Bspline surface_Ring}(b) for the obtained surface. If we choose the node vectors as $s_{ij}=s_i=s_{i-1}+1.5$, $i=1,\ldots,8$, and $t_{ij}=t_j=t_{j-1}+1.0$, $j=1,\ldots,20$, a $C^2$ moving B-spline surface with rounded feature lines is obtained; see Figure~\ref{Fig:moving Bspline surface_Ring}(c).

\subsection{Weighted moving B-spline surfaces}
Same as moving B-spline curves, moving B-spline surfaces can also have weights for the control points. We just illustrate here weighted moving B-spline surfaces defined on an open domain, closed weighted moving B-spline surfaces in one or two parameter directions can be defined similarly.

Suppose that $P_{ij}\in\mathbb{R}^d$, $i=0,1,\ldots,m$, $j=0,1,\ldots,n$, are a set of given points and $\omega_{ij}>0$, $i=0,1,\ldots,m$, $j=0,1,\ldots,n$, are the given weights. Assume the nodes corresponding to the points are $(s_{ij},t_{ij})$, $i=0,1,\ldots,m$, $j=0,1,\ldots,n$. A weighted moving B-spline surface is given by
\begin{equation}\label{Eqn:weighted_Moving_Bspline_surface}
\begin{array}{l}
P(s,t)=\frac{\sum_{i=0}^m\sum_{j=0}^n \omega_{ij}P_{ij} N_{k_1}(s-s_{ij})N_{k_2}(t-t_{ij})}{\sum_{i=0}^m\sum_{j=0}^n \omega_{ij}N_{k_1}(s-s_{ij})N_{k_2}(t-t_{ij})},
%(s,t)\in[s_1-\frac{k_1}{2},s_{m-1}+\frac{k_1}{2}]\times[t_1-\frac{k_2}{2},t_{n-1}+\frac{k_2}{2}].
\end{array}
\end{equation}
where the parameter domain is the same as that given by Equation (\ref{Eqn:Moving_Bspline_surface}).
Similar to the curve case, the weights will be used as additional tools for shape control.

\begin{figure}[htb]
  \centering
  \subfigure[]{\includegraphics[width=0.3\linewidth]{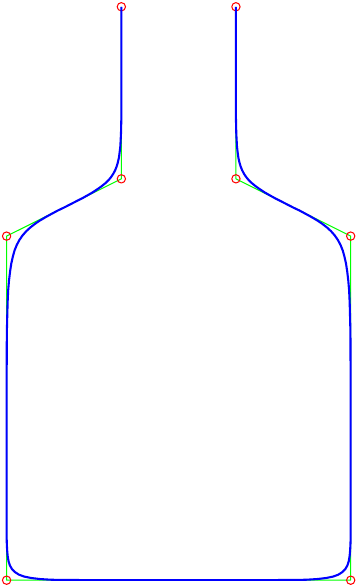}}\\
  \subfigure[]{\includegraphics[width=0.49\linewidth]{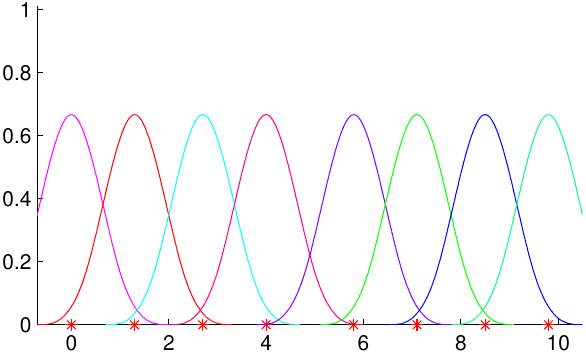}}
  \subfigure[]{\includegraphics[width=0.49\linewidth]{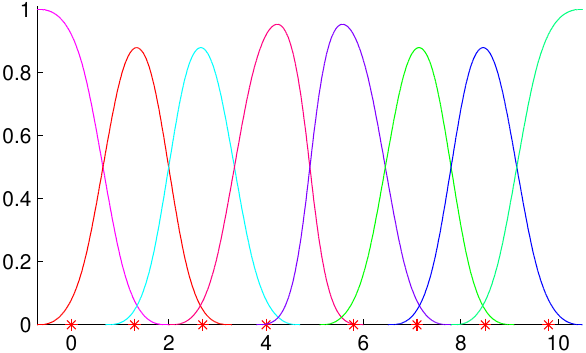}}
  \caption{An open moving B-spline curve of degree 3 with non-uniform nodes: (a) the moving B-spline curve; (b) the moving cubic B-spline bases; (c) the normalized moving B-spline bases. }
  \label{Fig:moving Bspline curves_bottle}
\end{figure}

\begin{figure}[htb]
  \centering
  \includegraphics[width=0.3\linewidth]{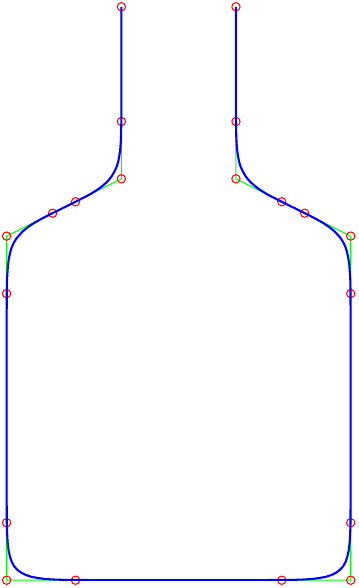}
  \caption{The bottle like shape is modeled by a uniform cubic B-spline curve.}
  \label{Fig:bottle_cubic}
\end{figure}

%%%%%%%%%%%%%%%%%%%%%%%%%%%%%%%%%%%%%%%%%%%%%%%%%%%%%%%%%%%%%%%%%%%%%%%%%%%%%%%%%%%%%%%%%%%
%%% Sectioin 4                                                                          %%%
%%%%%%%%%%%%%%%%%%%%%%%%%%%%%%%%%%%%%%%%%%%%%%%%%%%%%%%%%%%%%%%%%%%%%%%%%%%%%%%%%%%%%%%%%%%

\section{Examples}
\label{Sec:Examples}

In this section we present several interesting examples to show how to model curves and surfaces with salient features using moving B-splines. Comparisons with some existing modeling techniques will also be given or discussed.

\begin{figure}[htb]
  \centering
  \subfigure[]{\includegraphics[width=0.35\linewidth]{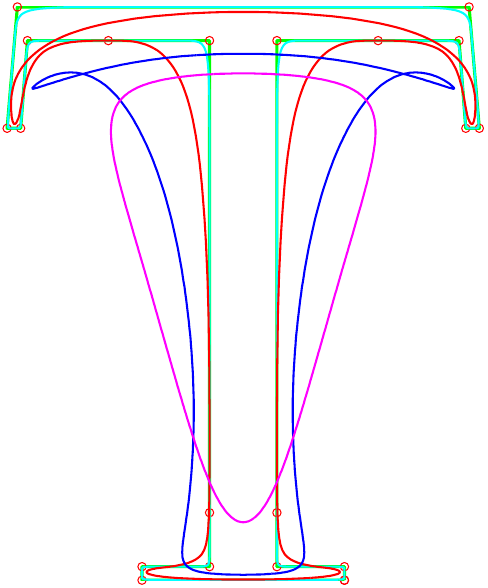}}
  \subfigure[]{\includegraphics[width=0.49\linewidth]{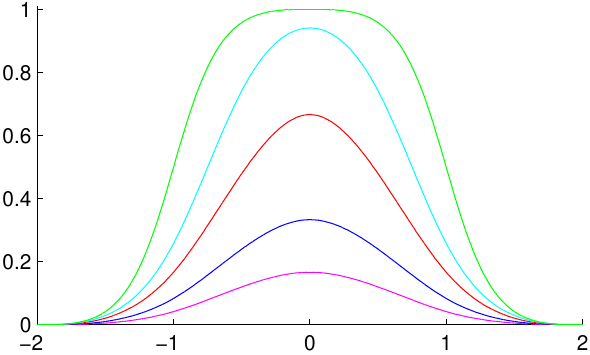}}
  \subfigure[]{\includegraphics[width=0.35\linewidth]{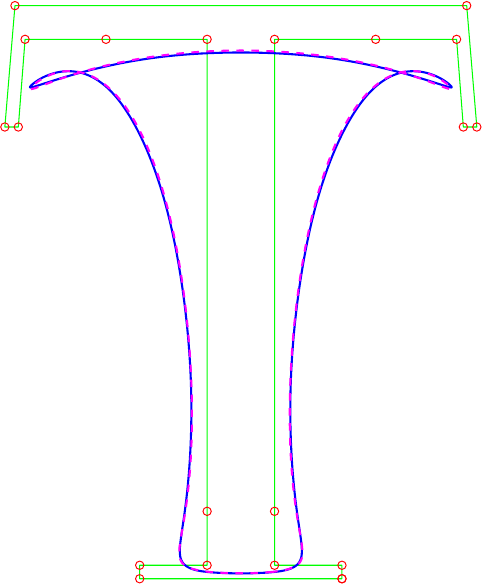}}
  \caption{(a) Uniform moving B-spline curves with various node intervals; (b) the normalized moving cubic B-spline bases with node intervals (from top to bottom) 2, 1.5, 1.0, 0.5, 0.25, respectively; (c) uniform moving B-spline curve of degree 3 (solid) with knot interval 0.5 and a uniform B-spline curve of degree 15 (dashed). }
  \label{Fig:moving Bspline multicurves}
\end{figure}

\begin{figure}[htb]
  \centering
  \subfigure[]{\includegraphics[width=0.35\linewidth]{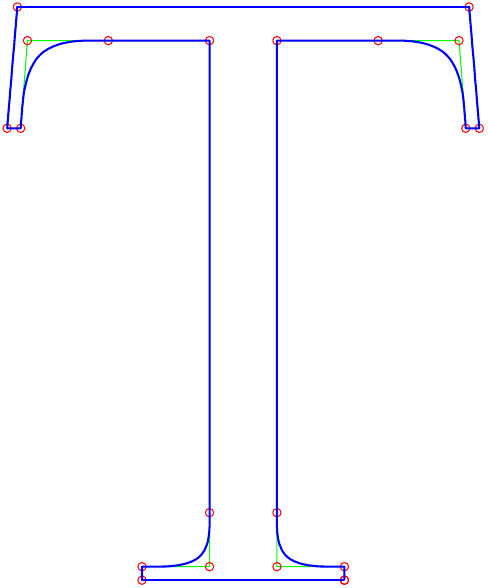}}\\
  \subfigure[]{\includegraphics[width=0.49\linewidth]{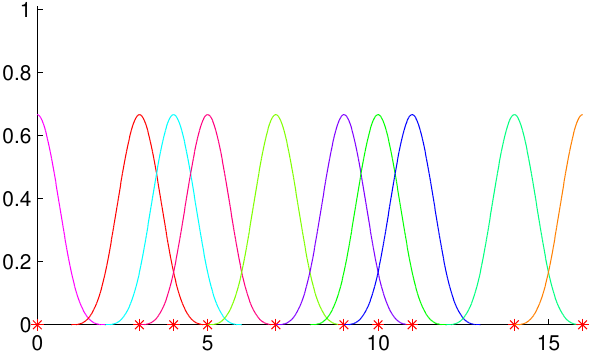}}
  \subfigure[]{\includegraphics[width=0.49\linewidth]{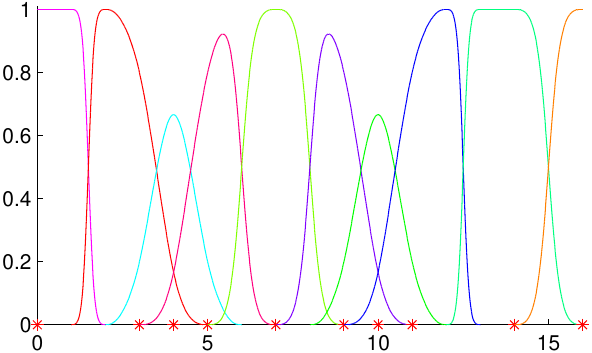}}
  \caption{A moving B-spline curve of degree 3 with non-uniform nodes: (a) the moving B-spline curve; (b) the moving cubic B-spline bases; (c) the normalized moving B-spline bases. }
  \label{Fig:moving Bspline curves_T}
\end{figure}

First, we model a bottle like shape using moving B-splines. Given a control polygon consisting of 8 control points together with nodes (0, 1.3, 2.7, 4.0, 5.8, 7.1, 8.5, 9.8), a cubic moving B-spline curve is obtained; see Figure~\ref{Fig:moving Bspline curves_bottle}(a) for the obtained curve. The moving cubic B-spline bases and the normalized bases are given in Figures~\ref{Fig:moving Bspline curves_bottle}(b) and \ref{Fig:moving Bspline curves_bottle}(c), respectively. Because there exist parameter intervals on which only two basis functions do not vanish, the curve has local straight line segments on every edge. If a similar shape is modeled by conventional cubic B-spline curve, more knots will be inserted and more control points will be used. See Figure~\ref{Fig:bottle_cubic} for the bottle like shape modeled by a uniform cubic B-spline curve. Note that the bottle like shape can also be modelled by MD-spline with a small number of control points. The difference lies at that the moving B-spline curve can have a consistent high order of smoothness along the curve while the smoothness orders at all knots of a MD-spline may be different.

\begin{figure}[htb]
  \centering
  \subfigure[]{\includegraphics[width=0.49\linewidth]{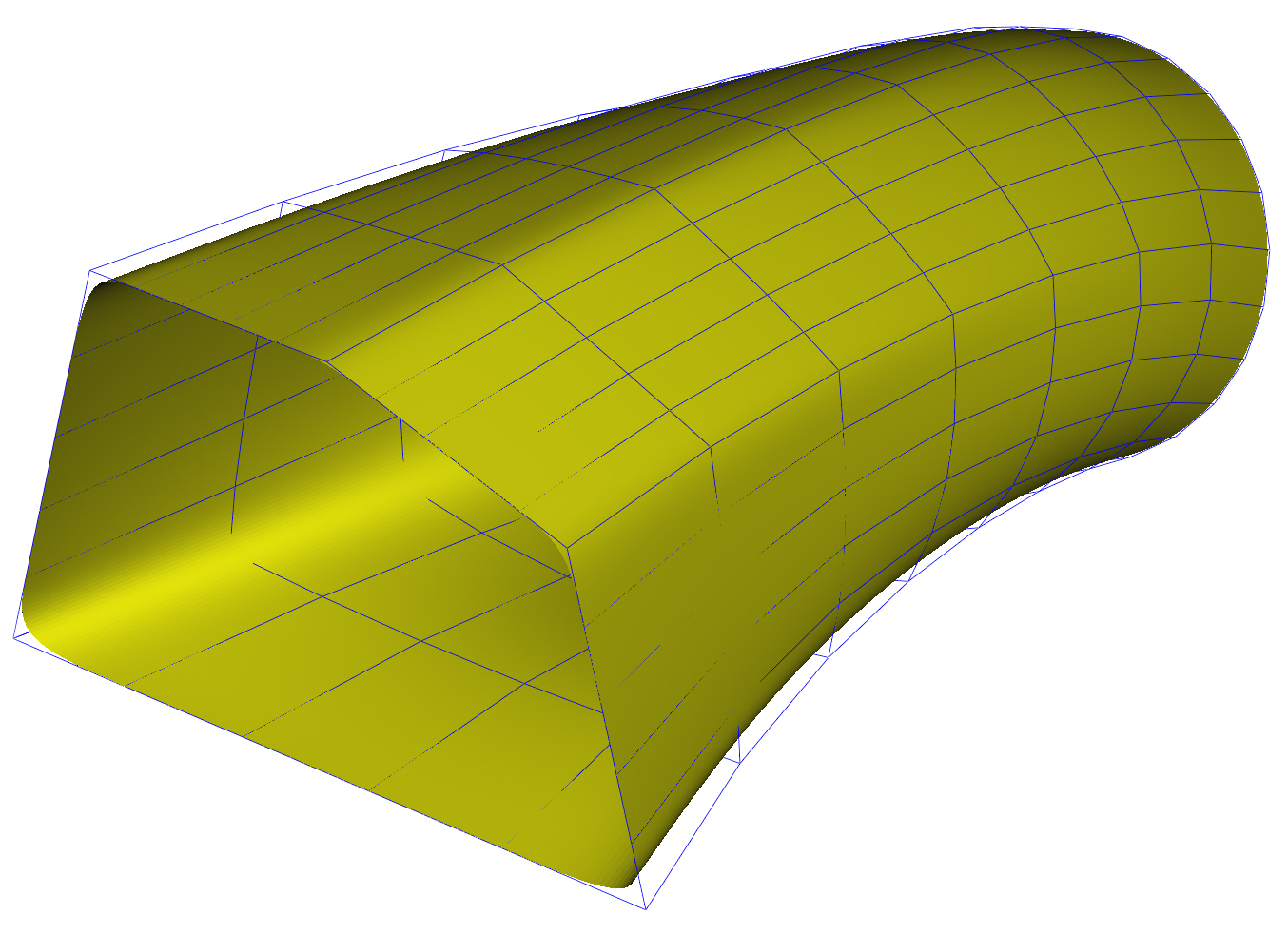}}
  \subfigure[]{\includegraphics[width=0.49\linewidth]{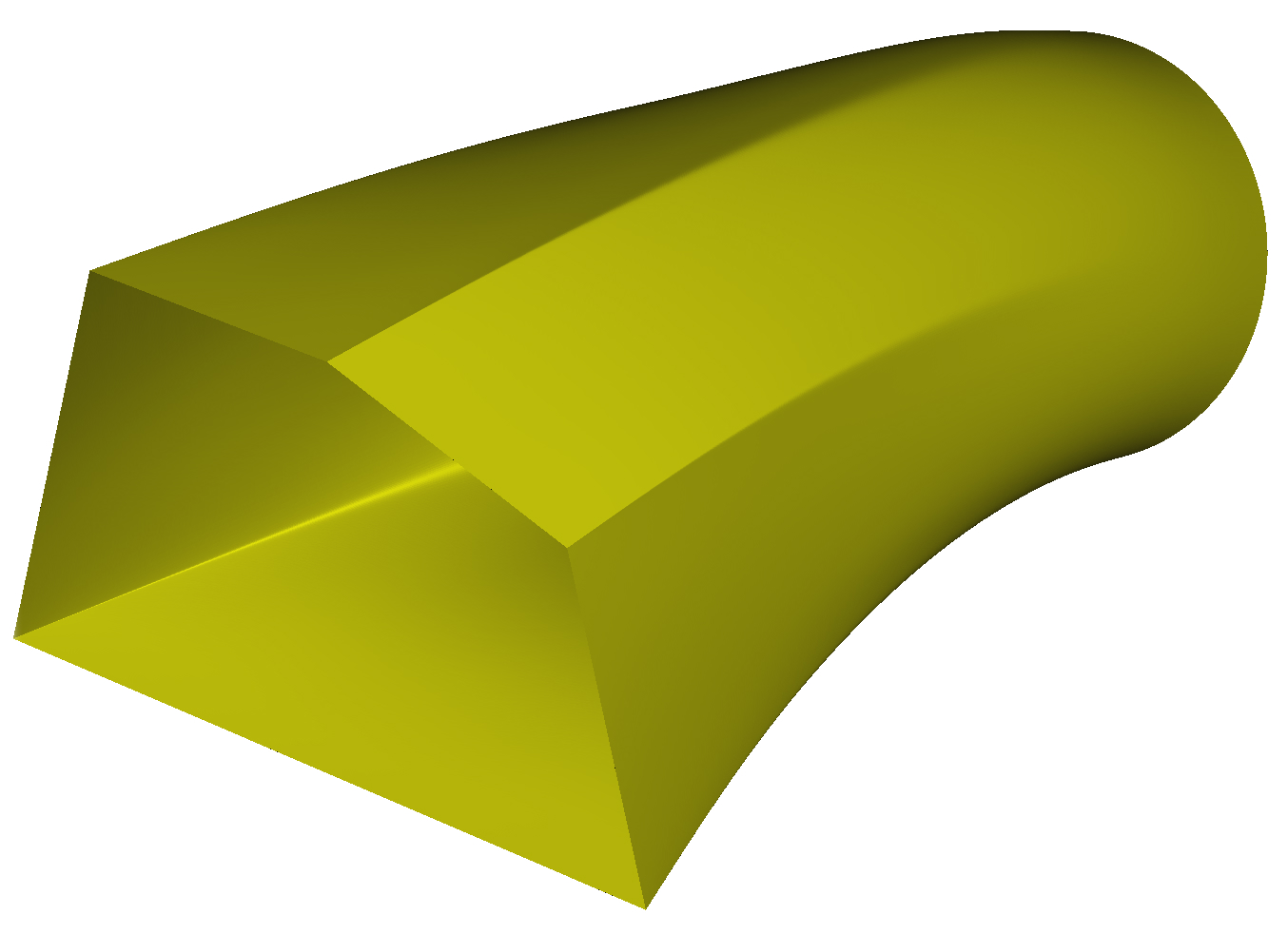}}
  \caption{Modeling a tunnel shape by moving B-splines: (a) the input control mesh and the bicubic B-spline surface; (b) the surface constructed by bicubic moving B-splines. }
  \label{Fig:moving bicubic_tunnel}
\end{figure}

Second, we present an example to show how a polygon can be reproduced or smoothed using a moving cubic B-spline. Figure~\ref{Fig:moving Bspline multicurves}(a) illustrates a closed control polygon consisting of 20 control points. Suppose that the nodes are uniformly sampled as $t_0=0$, $t_i=t_{i-1}+a$, $i=1,2,\ldots,20$. Five cubic moving B-spline curves are constructed by choosing the node intervals as $a=2.0,1.5,1.0,0.5,0.25$, respectively. Figure~\ref{Fig:moving Bspline multicurves}(b) illustrates the normalized moving B-spline bases corresponding to various node intervals. From the figure we know that the original control polygon is reproduced exactly when the node intervals are chosen as $a=2.0$ for the moving cubic B-spline and the control polygon is smoothed when the node intervals are reduced. In particular, moving B-spline curves with a fixed degree but various node intervals can approximate conventional B-spline curves with different degrees. Figure~\ref{Fig:moving Bspline multicurves}(c) illustrates that the cubic moving B-spline curve with the same node interval $a=0.5$ can approximate a uniform B-spline curve of degree 15 very well.
Assume the bottom right corner is $P_0$ and the control points are numbered counter-clockwise. By choosing the nodes as (0,3,4,5,7,9,10,11,14,16,18,20,23,24,25,27,29,30,32,34,36), a moving cubic B-spline curve with non-uniform nodes is obtained; see Figure~\ref{Fig:moving Bspline curves_T} for the obtained moving B-spline curve and the partial bases. Though other techniques such as multi-degree B-splines can also be employed to model similar shape, the moving B-spline curves own the advantages of simplicity and universality.

Third, we model a tunnel shape with moving bicubic B-splines. Figure~\ref{Fig:moving bicubic_tunnel}(a) illustrates the input control mesh that is consisting of $8\times25$ control points. The control points on one boundary lie on a pentagon while the control points on the other one lie on a circle. When a standard bicubic B-spline surface is constructed from the control mesh, all sharp corners on the control mesh have been smoothed out. To model a tunnel shape surface that has a piecewise linear boundary on one side and a circular boundary on the other side, we choose nodes adaptively for the construction of a bicubic moving B-spline surface. Let $s_i=i$, $i=0,1,\ldots,7$. The node parameters $s_{ij}$s for the control points $P_{ij}$ are chosen as $s_{ij}=s_i$, $i=0,1,\ldots,7$, $j=0,1,\ldots,25$. Let $d_i=(\frac{i}{7})^2$, $i=0,1,\ldots,7$. The node parameters $t_{ij}$ for the control points are given by
\begin{eqnarray*}
    t_{ij}=
    \left\{\begin{array}{ll}
    j &                 \textrm{if} \ j \ \texttt{mod}\ 5=0  \\
    j+d_i   &           \textrm{if} \ j \ \texttt{mod}\ 5=1  \\
    j+\frac{d_i}{3} &   \textrm{if} \ j \ \texttt{mod}\ 5=2  \\
    j-\frac{d_i}{3} &   \textrm{if} \ j \ \texttt{mod}\ 5=3  \\
    j-d_i   &           \textrm{if} \ j \ \texttt{mod}\ 5=4
    \end{array}\right.
\end{eqnarray*}
When all nodes have been chosen, a surface that is open in $s$ direction but closed in $t$ direction is obtained by Equation~(\ref{Eqn:Moving_Bspline_surface_close_v}).
Based on the choice of the nodes we know that $t_{ij}=j$, $j=0,1,\ldots,25$ when $i=0$. Thus the surface has a smooth boundary at $s=s_1-2.0$. If $i=7$, the node parameters satisfy $t_{7,j-1}=j-2$, $t_{7,j}=j$, and $t_{7,j+1}=j+2$ when $j \ \texttt{mod}\ 5=0$. The obtained surface interpolates all sharp corners $P_{7,j}$, $j=0,5,\ldots,20$, and the surface boundary is a pentagon as defined by the control mesh.

%%%%%%%%%%%%%%%%%%%%%%%%%%%%%%%%%%%%%%%%%%%%%%%%%%%%%%%%%%%%%%%%%%%%%%%%%%%%%%%%%%%%%%%%%%%
%%% Sectioin 5                                                                          %%%
%%%%%%%%%%%%%%%%%%%%%%%%%%%%%%%%%%%%%%%%%%%%%%%%%%%%%%%%%%%%%%%%%%%%%%%%%%%%%%%%%%%%%%%%%%%

\section{Conclusions and discussions}
\label{Sec:Conclude}
This paper has proposed a novel technique of constructing rational curves and surfaces by uniform B-splines. Given a control polygon or a control mesh together with node ordinates corresponding to all control points, a rational curve or surface is obtained by normalized moving B-splines or normalized weighted moving B-splines centered at the given nodes. The obtained moving B-spline curves can have sharp or rounded corners, partial or full straight edges while the obtained surfaces can have sharp or rounded vertices, sharp or smoothed edges, etc. Users can design curves or surfaces with prescribed features just by choosing proper nodes for the control points. Compared with conventional NURBS curves and surfaces that represent the similar shapes, the curves and surfaces constructed by the proposed technique have less control points. The proposed technique also benefits much from simplicity and easy implementation for practical applications.

At present, we have constructed moving B-spline curves and surfaces with ordered control points. The nodes are also assumed monotone for curves and the union of all parameter domains of moving tensor product B-splines form a rectangle. Since curves and surfaces constructed with moving B-splines are point based, it is possible to model discontinuous curves with unorganized control points or surfaces with more complex parameter domains by the proposed technique. From the second example illustrated Section~\ref{Sec:Examples} we know that the moving B-spline curves in terms of the parameter $t$ and the uniform node interval $a$ can form a continuous domain bounded by a closed polygon. A bivariate planar domain bounded by a closed polygon may find applications such as free-form deformation or isogeometric analysis, etc.

\section*{Acknowledgment}

This work was supported by the National Natural Science Foundation of China under grant No. 12171429.

\bibliographystyle{elsarticle-harv}
\bibliography{movingBsplineCurvesSurfaces}




%% \bibitem must have one of the following forms:
%%   \bibitem[Jones et al.(1990)]{key}...
%%   \bibitem[Jones et al.(1990)Jones, Baker, and Williams]{key}...
%%   \bibitem[Jones et al., 1990]{key}...
%%   \bibitem[\protect\citepauthoryear{Jones, Baker, and Williams}{Jones
%%       et al.}{1990}]{key}...
%%   \bibitem[\protect\citepauthoryear{Jones et al.}{1990}]{key}...
%%   \bibitem[\protect\astroncite{Jones et al.}{1990}]{key}...
%%   \bibitem[\protect\citepname{Jones et al., }1990]{key}...
%%   \harvarditem[Jones et al.]{Jones, Baker, and Williams}{1990}{key}...
%%


% \end{thebibliography}

\end{document}